\documentclass[reqno,12pt]{amsart}
\usepackage{amscd,amssymb}

\theoremstyle{plain}
\newtheorem{Thm}[subsection]{Theorem}
\newtheorem{Cor}[subsection]{Corollary}
\newtheorem{Lem}[subsection]{Lemma}
\newtheorem{Prop}[subsection]{Proposition}
\newtheorem{Conj}[subsection]{Conjecture}

\theoremstyle{definition}
\newtheorem{Def}[subsection]{Definition}

\theoremstyle{remark}

\newtheorem{Rem}[subsection]{Remark}

\newtheorem{conj}{Conjecture}

\numberwithin{equation}{section}
\renewcommand{\rm}{\normalshape}

\newif\ifShowLabels
\ShowLabelstrue
\newdimen\theight
\def\TeXref#1{%
        \leavevmode\vadjust{\setbox0=\hbox{{\tt
                \quad\quad  {\small \rm #1}}}%
        \theight=\ht0
        \advance\theight by \lineskip
        \kern -\theight \vbox to
        \theight{\rightline{\rlap{\box0}}%
        \vss}%
        }}%

\ShowLabelsfalse% comment this out if labels should be printed

\renewcommand{\sec}[2]{\section{#2}\label{S:#1}%
        \ifShowLabels \TeXref{{S:#1}} \fi}
\newcommand{\ssec}[2]{\subsection{#2}\label{SS:#1}%
        \ifShowLabels \TeXref{{SS:#1}} \fi}

\newcommand{\refs}[1]{Section ~\ref{S:#1}}
\newcommand{\refss}[1]{Section ~\ref{SS:#1}}

\newcommand{\reft}[1]{Theorem ~\ref{T:#1}}
\newcommand{\refl}[1]{Lemma ~\ref{L:#1}}
\newcommand{\refp}[1]{Proposition ~\ref{P:#1}}
\newcommand{\refc}[1]{Corollary ~\ref{C:#1}}

\newcommand{\refe}[1]{\eqref{E:#1}}
\newcommand{\refco}[1]{Conjecture ~\ref{Co:#1}}

\newenvironment{thm}[1]%
        { \begin{Thm} \label{T:#1}  \ifShowLabels \TeXref{T:#1} \fi }%
        { \end{Thm} }

\renewcommand{\th}[1]{\begin{thm}{#1} \sl }
\renewcommand{\eth}{\end{thm} }

\newenvironment{lemma}[1]%
        { \begin{Lem} \label{L:#1}  \ifShowLabels \TeXref{L:#1} \fi }%
        { \end{Lem} }
\newcommand{\lem}[1]{\begin{lemma}{#1} \sl}
\newcommand{\elem}{\end{lemma}}

\newenvironment{propos}[1]%
        { \begin{Prop} \label{P:#1}  \ifShowLabels \TeXref{P:#1} \fi }%
        { \end{Prop} }
\newcommand{\prop}[1]{\begin{propos}{#1}\sl }
\newcommand{\eprop}{\end{propos}}

\newenvironment{corol}[1]%
        { \begin{Cor} \label{C:#1}  \ifShowLabels \TeXref{C:#1} \fi }%
        { \end{Cor} }
\newcommand{\cor}[1]{\begin{corol}{#1} \sl }
\newcommand{\ecor}{\end{corol}}

\newenvironment{defeni}[1]%
        { \begin{Def} \label{D:#1}  \ifShowLabels \TeXref{D:#1} \fi }%
        { \end{Def} }
\newcommand{\defe}[1]{\begin{defeni}{#1} \sl }
\newcommand{\edefe}{\end{defeni}}

\newenvironment{remark}[1]%
        { \begin{Rem} \label{R:#1}  \ifShowLabels \TeXref{R:#1} \fi }%
        { \end{Rem} }
\newcommand{\rem}[1]{\begin{remark}{#1}}
\newcommand{\erem}{\end{remark}}

\newenvironment{conjec}[1]%
        { \begin{Conj} \label{Co:#1}  \ifShowLabels \TeXref{Co:#1} \fi }%
        { \end{Conj} }
\renewcommand{\conj}[1]{\begin{conjec}{#1} \sl }
\newcommand{\econj}{\end{conjec}}

\newcommand{\eq}[1]%
        { \ifShowLabels \TeXref{E:#1} \fi
           \begin{equation} \label{E:#1} }
\newcommand{\eeq}{ \end{equation} }

\newcommand{\prf}{ \begin{proof} }
\newcommand{\epr}{ \end{proof} }

\newcommand\alp{\alpha}         

\newcommand\gam{\gamma}         
\newcommand\del{\delta}

\newcommand\kap{\kappa}
\newcommand\lam{\lambda}                
\newcommand\sig{\sigma}

\newcommand\calD{{\mathcal{D}}}

\newcommand\calF{{\mathcal{F}}}
\newcommand\calG{{\mathcal{G}}}
\newcommand\calH{{\mathcal{H}}}

\newcommand\calK{{\mathcal{K}}}
\newcommand\calL{{\mathcal{L}}}

\newcommand\calN{{\mathcal{N}}}

\newcommand\bfa{{\mathbf a}}            
            \newcommand\bfB{{\mathbf B}}
            \newcommand\bfC{{\mathbf C}}

\newcommand\bff{{\mathbf f}}            \newcommand\bfF{{\mathbf F}}
            \newcommand\bfG{{\mathbf G}}
            \newcommand\bfH{{\mathbf H}}
\newcommand\bfi{{\mathbf i}}            
\newcommand\bfj{{\mathbf j}}            
            
            \newcommand\bfL{{\mathbf L}}
\newcommand\bfm{{\mathbf m}}            \newcommand\bfM{{\mathbf M}}

\newcommand\bfp{{\mathbf p}}            \newcommand\bfP{{\mathbf P}}
\newcommand\bfq{{\mathbf q}}            
\newcommand\bfr{{\mathbf r}}            
\newcommand\bfs{{\mathbf s}}            \newcommand\bfS{{\mathbf S}}
            \newcommand\bfT{{\mathbf T}}
            \newcommand\bfU{{\mathbf U}}

            \newcommand\bfX{{\mathbf X}}
            
            \newcommand\bfZ{{\mathbf Z}}

\newcommand\QQ{\mathbb{Q}}

\renewcommand\AA{\mathbb{A}}

\newcommand\DD{\mathbb{D}}
\newcommand\FF{\mathbb{F}}
\newcommand\GG{\mathbb{G}}

\newcommand\ZZ{\mathbb{Z}}

 \newcommand\grb{{\mathfrak{b}}}

 \newcommand\grg{{\mathfrak{g}}}

 \newcommand\grt{{\mathfrak{t}}}

\newcommand\sdp{\times \hskip -0.3em {\raise 0.3ex
\hbox{$\scriptscriptstyle |$}}} % semidirect product

\newcommand\End{\operatorname{End\,}}
\newcommand\Ext{\operatorname{Ext}}

\newcommand\gl{{\bf gl}}

\newcommand\Hom{\operatorname {Hom}}
\newcommand\RHom{\operatorname {RHom}}

\newcommand\Id{\operatorname{Id}}
\newcommand\im{\operatorname {im}}

\newcommand\Ind{\operatorname{Ind}}

\newcommand\Res{\operatorname{Res}}

\newcommand{\sign}{\operatorname{sign}}

\newcommand\supp{\operatorname{supp}}

\newcommand\tr{\operatorname{tr}}
\newcommand\Tr{\operatorname{Tr}}

\newcommand\oF{{\overline{F}}}

\newcommand\x{\times}
\newcommand\ten{\otimes}

\renewcommand{\>}{\rangle}
\newcommand{\<}{\langle}

\newcommand{\ra}{\rangle}
\newcommand{\la}{\langle}

\newcommand\qlb{{\overline \QQ}_l}
\newcommand\fr{\text{Fr}}
\newcommand\frw{\fr_w}
\newcommand\gln{{\bf GL}(n)}
\newcommand\wt{\widetilde}
\newcommand\chl{\calK_{\calL}}
\newcommand\ch{\text{ch}}
\newcommand\fq{\FF_q}

\newcommand{\bPhi}{\Phi\hskip-7.5pt \Phi}

\renewcommand{\gl}{{\bf GL}}

\renewcommand{\Id}{\text{Id}}
\newcommand{\Irr}{\operatorname{Irr}}

\newcommand\tbG{{\widetilde \bfG}}

\newcommand\td{\bfT^{\vee}}

\newcommand\phrog{\bPhi_{\bfG,\rho,\psi}}
\newcommand\phrot{\bPhi_{\bfT,\rho,\psi}}
\newcommand\kg{\text{K}_G}
\newcommand\obfT{{\overline \bfT}}
\newcommand\Spr{{\bf Spr}}

\newcommand\avgb{{\bf Av_{G/B}}}
\newcommand\tilInd{{\widetilde \Ind}}
\newcommand\tilRes{{\widetilde \Res}}

\newcommand\qlbs{(\qlb[1](\frac{1}{2}))}
\newcommand\tilgrg{{\widetilde \grg}}
\begin{document}

\title{$\gam$-sheaves on reductive groups}
\author{Alexander Braverman and David Kazhdan}
\begin{abstract}
Let $\bfG$ be a reductive group over a finite field $F=\FF_q$. Fix a non-trivial additive character
$\psi:F\to\qlb^\x$. In \cite{BrKa}
we introduced certain $\gam$-functions $\gam_{G,\rho,\psi}$ on the set $\Irr(G)$ of irreducible
representations of the finite group $G=\bfG(F)$.

As usual every function $\gam_{G,\rho,\psi}$ on $\Irr(G)$ gives rise to an
Ad$G$-equivariant function $\Phi_{G,\rho,\psi}$ on $G$. The purpose of this paper
is to construct an irreducible perverse sheaf $\phrog$ on $\bfG$ such that
the function $\Phi_{G,\rho,\psi}$ is obtained conjecturally by taking traces of Frobenius
morphism in the stalks of $\phrog$. In order to do this we need to assume that
$\rho$ satisfies certain technical condition (we call $\rho$ {\it good} if that condition is
satisfied). We prove this conjecture for $\bfG=\gl(n)$ and for $\bfG$ of
semi-simple rank one.

We also prove the above conjecture assuming that certain cohomology vanishing for the sheaf $\phrog$
holds (we show that this is the case  for
$\bfG$ of semisimple rank 1).  Assuming this vanishing
we show that if both $\rho_1$ and $\rho_2$ are good then
$\bPhi_{\bfG,\rho_1,\psi}\star \bPhi_{\bfG,\rho_2,\psi}\simeq \bPhi_{\bfG,\rho_1\oplus\rho_2,\psi}$.
We also compute the convolution of $\phrog$  with the majority of Lusztig's character sheaves.

We conjecture that the functor of convolution with $\phrog$ is exact in the perverse $t$-structure.
\end{abstract}
\maketitle

%-------------------------------------------------------------------------------------------------------------
\sec{int}{Introduction}
%---------------------------------------------------------------------------------------
\ssec{somnot}{Some notations}
In what follows we denote by $F=\fq$ the finite field with $q$ elements, $\oF$ -- its algebraic closure.
Choose a prime number $\ell$ which is prime to $q$. Let also $\psi:F\to
\qlb^\x$ denote a non-trivial additive character of $F$ with values in
$\qlb^\x$. We shall denote algebraic varieties over $F$ by boldface letters
(e.g. $\bfG$, $\bfX$ etc).
The corresponding ordinary letters (e.g. $G$, $T$ etc) will denote the
corresponding sets of $F$-points.

For a finite group $G$ we denote by $\Irr(G)$ the set of isomorphism classes
of irreducible representations of $G$.

In what follows we choose a square root $q^{1/2}$ of $q$.

For an algebraic variety $\bfX$ over $\oF$ we shall denote by
$\calD(\bfX)$ the bounded derived category of $\ell$-adic sheaves on
$\bfX$.

If $X$ is defined over $F$ we let $\fr:\bfX\to\bfX$ denote the geometric
Frobenius morphism. We say that an object $\calF\in \calD(\bfX)$ is endowed with
a Weil structure if we are given an isomorphism
$\fr^*\calF\simeq\calF$. To any Weil sheaf on $\bfX$ we associate a function $\chi(\calF)$
on $X=\bfX(F)$ in the  following way. Let $x\in X$ and let $\calF_x$ denote
the fiber of $\calF$ at $x$. This is a complex of $\ell$-adic vector spaces.
Since $x$ is fixed by $\fr$,
the Weil structure on $\calF$ gives rise to an automorphism of $\calF_x$ which
(by abuse of language) we shall also denote by $\fr$. Thus we set
\eq{}
\chi(\calF)(x)=\sum_i (-1)^i \Tr(\fr,H^i(\calF_x))
\end{equation}

Let $\calF$ be a Weil sheaf. For a half integer $n$ we denote
by $\calF(n)$ the Tate twist of $\calF$ (corresponding to the
chosen $q^{1/2}$). Thus $\chi(\calF(n))=\chi(\calF)q^{-n}$.

For an algebraic group $\bfG$ we define two convolution functors
$(\calF,\calG)\mapsto \calF\star\calG$ and $(\calF,\calG)\mapsto \calF*\calG$
going from $\calD(\bfG)\x \calD(\bfG)$ to $\calD(\bfG)$ in the following way.
Let $\bfm:\bfG\x\bfG\to\bfG$ denote the multiplication map. Then
\eq{}
\calF\star\calG=\bfm_!(\calF\boxtimes\calG)\quad\text{and}\quad
\calF*\calG=\bfm_*(\calF\boxtimes\calG).
\end{equation}
%---------------------------------------------------------------------------------------------
\ssec{gln-int}{$\gam$-functions for $\gl(n)$}
 Let $(\pi,V)$ be an irreducible
representation of $G=\gl(n,F)$. Consider the operator
\eq{fing}
\sum\limits_{g\in G}\psi(\tr(g))\pi(g)(-1)^n q^{-n^2/2}
\in \End_G V.
\end{equation}
By Schur's lemma this operator takes the form
$\gam_\psi(\pi)\cdot \Id_V$ where $\gam_\psi(\pi)\in \qlb$.

The number $\gam_\psi(\pi)$ is called the gamma-function
of the representation $\pi$. One can "explicitly" compute
$\gam_\psi(\pi)$ in the following way.

Let $W\simeq S_n$ denote the Weyl group of $\gl(n)$. Following Deligne and
Lusztig (cf. \cite{dl}) we can associate to every $w$ a maximal torus
$\bfT_w\subset \gl(n)$ (defined uniquely up to $G=\gl(n,F)$-conjugacy).

Fix $w\in W$. For a character $\theta:T_w\to \qlb^{\x}$ we set
\eq{}
\gam_{\psi,w}(\theta)=(-1)^n q^{-n/2}\sum\limits_{t\in
T_w}\psi(\tr(t))\theta(t)\in\qlb.
\end{equation}
\medskip
{\bf Example.} Assume that $w\in S_n$ is a cycle of length $n$.
Then $T_w\simeq E^{\x}$ where $E$ is the (unique up to
isomorphism) extension of $F$ of degree $n$. In this case
$\gam_w(\theta)=(-q^{1/2})^{n-1}\gam_{E,\psi}(\theta)$
for any character $\theta$ of
$E^{\x}$, where by $\gam_{E,\psi}(\theta)$ we denote the $\gam$-function
defined as in \refe{fing} for the group $\gl(1,E)\simeq E^{\x}$.

Recall that in \cite{dl} Deligne and Lusztig have associated to
$\theta$ a virtual representation $R_{\theta,w}$ of $G$ and they have
proved that every $\pi\in \Irr(G)$ is an irreducible constituent of some
(in general non-unique) $R_{\theta,w}$.

The following result is proven in \refs{gln}.
%------------------------------------------------------
\th{gln-func}Assume that an irreducible representation
$(\pi,V)$ appears in $R_{\theta,w}$ for some
$w$ and $\theta$ as above. Then
\eq{}
\gam_\psi(\pi)=\gam_{w,\psi}(\theta)
\end{equation}
In particular, $\gam_\psi(\pi)=\gam_\psi(\pi')$ if $\pi$ and $\pi'$
appear in the same virtual representation $R_{\theta,w}$.
\eth
%-------------------------------------------------------------------------
\ssec{int-arb-group}{The case of arbitrary group}
Let now $\bfG$ be any connected split
reductive group over $F$, $G=\bfG(F)$.
Let $\bfT$ be the Cartan group of $\bfG$. Let also $\td$ denote the dual torus
to $\bfT$ over $\qlb$. The Weyl group $W$ acts naturally on $\td$.

Assume that we are given an $n$-dimensional representation $\rho:\td\to\gl(n,\qlb)$ of
$\td$ such that for every $w\in W$ the composition $\rho\circ w$ is isomorphic to $\rho$. In other words,
$\rho$ is given by a collection $\lam_1,...,\lam_n$ of characters of $\td$ which is invariant
under the action of $W$.

Let $\td_\rho=\GG_{m,\qlb}^n$. Then we get a natural map
$\bfp_\rho^{\vee}:\td\to\td_\rho$ sending every $t$ to
$(\lam_1(t),...,\lam_n(t))$.

Let now $\bfT_\rho\simeq \GG_{m,F}^n$ denote the dual torus to
$\td_\rho$ over $F$ and let $\bfp_\rho:\bfT_\rho\to \bfT$ denote
the map which is dual to $\bfp_\rho^{\vee}$. Explicitly one has
\eq{}
\bfp_\rho(x_1,...,x_n)=\lam_1(x_1)...\lam_n(x_n) .
\end{equation}

Let $W_\rho\simeq S_n$ denote the Weyl group of $\gl(n)$.

Let now $\pi$ be an irreducible representation of $G$. Assume that
$\pi$ appears in some $R_{\theta,w}$ for some $\theta:T_w\to
\qlb^{\x}$. Let $w'$ be any lift of $w$ to $W_\rho$. Then
$\bfp_\rho$ induces an $F$-rational map
$\bfp_{w'}:\bfT_{\rho,w'}\to\bfT_w$, hence a homomorphism
$p_{w'}:T_{\rho,w'}\to T_w$. Let $\theta'=p_{w'}^*(\theta)$. Define
    \eq{}
    \gam_{G,\rho,\psi}(\pi):=\gam_\psi(\pi'),
    \end{equation}
where $\pi'$ is any irreducible representation of
$G_\rho$ which appears in $R_{\theta',w'}$. By
\reft{gln-func} one has
\eq{}
\gam_{G,\rho,\psi}(\pi)=\gam_{w',\psi}(\theta').
\end{equation}
\lem{}
The definition of $\gam_{G,\rho,\psi}(\pi)$ does not depend on the
choice of $w'$.
\elem
\prf Let $w''$ be another lift of $w$ to $W_\rho$ and let $\theta''$
be the corresponding character of $T_{w''}$.
Then  it follows from \cite{dl} that
\eq{}
\la R_{\theta',w'},R_{\theta'',w''}\ra\neq 0.
\end{equation}
Therefore, our lemma follows from \reft{gln-func}.
\epr
Sometimes we shall write $\gam_{w,\rho,\psi}(\theta)$ instead of $\gam_{G,\rho,\psi}(\pi)$.
%--------------------------------------------------------------------
Let now $\Phi_{G,\rho,\psi}$ denote the unique central function on $G$
such that for every irreducible representation $(\pi,V)$ of $G$
one has
\eq{}
\sum_{g\in
G}\Phi_{G,\rho,\psi}(g)\pi(g)=\gam_{G,\rho,\psi}(\pi)\cdot\Id_V.
\end{equation}
We would like to compute this function explicitly using geometry.
More precisely, we are going to do the following.

We say that the representation $\rho$ is good
if there exists a character $\sig:\bfG\to\GG_m$ such that
for every weight $\lam_i$ of $\rho$ as above one
has $\la \lam_i,\sig\ra>0$.
For any good representation $\rho$ we are going to construct an irreducible
perverse sheaf $\phrog$ on $\bfG$ endowed with a Weil structure.

\noindent
{\it Remark.} The condition of being ``good'' is not very restrictive: if one starts with
arbitrary $\bfG$ and $\rho$ one cam always make $\rho$ good by
passing to $\bfG'=\bfG\x \GG_m$ and taking $\rho'=\rho\ten \text{St}$
where St denotes the standard one-dimensional representation of $\GG_m$.

One of our main results is
%-----------------------------------------------------------------------
\th{}Assume that the semi-simple rank of $\bfG$ is $\leq 1$ or that
$\bfG=\gl(n)$. Then
\eq{}
\chi(\phrog) =\Phi_{G,\rho,\psi}.
\end{equation}
\eth
%-----------------------------------------------------------
When $\rho$ is sufficiently generic (i.e. when the cocharacters $\lam_1,...,\lam_n$
span a lattice of rank equal to $\dim \bfT$)) the sheaf $\phrog$ as above is explicitly
constructed on the set $\bfG_r$
of regular elements in $\bfG$ and on the whole of $\bfG$ it is obtained
by means of the Goresky-MacPherson extension.

When the semi-simple rank of $\bfG$  is $\leq 1$
we can also show that the the functor of convolution with $\phrog$ enjoys
some nice properties. In \reft{main} we compute the convolution of $\phrog$ with Lusztig's
character sheaves. In particular, we show that $\star$ and $*$ convolutions
in this case coincide. Also, we prove the following result (assuming again that the semi-simple
rank of $\bfG$ is $\leq 1$):
\th{}
Assume that $\rho_1$ and $\rho_2$ are good with respect to the same
character $\sig$ of $\bfG$. Then
\eq{}
\bPhi_{\bfG,\rho_1\oplus\rho_2,\psi}\simeq
\bPhi_{\bfG,\rho_1,\psi}\star\bPhi_{\bfG,\rho_2,\psi}\simeq
\bPhi_{\bfG,\rho_1,\psi}*\bPhi_{\bfG,\rho_2,\psi}
\end{equation}
\eth

We conjecture that the above theorems hold for general $\bfG$ but we don't know how to prove this
(however, in \refs{main} we deduce these results from certain conjectural cohomology vanishing).
We also believe in the following
%----------------------------------------------------------------
\conj{}
The functors $\calF\mapsto\calF\star\phrog$ and $\calF\mapsto\calF*\phrog$ are
exact in the perverse $t$-structure.
\econj
For example, when $\bfG=\gl(n)$ and $\rho$ is the standard representation this conjecture
follows from the corresponding property of the Fourier-Deligne transform.
%------------------------------------------------------------------------------------
\ssec{ack}{Acknowlegdements}We are grateful to D.~Gaitsgory and R.~Bezrukavnikov, N.~Katz and
F.~Loeser
for very helpful discussions on the subject.
%--------------------------------------------------------------------------------------------------------------
\sec{}{Induction and restriction functors}
The purpose of this section is to collect some facts about Lusztig's induction and restriction
functors which will be used later.
%-------------------------------------------------------------------------------------------------
\ssec{restriction}{Restriction}Let $\bfP$ be a parabolic subgroup of $\bfG$ and let $\bfM$ be the
corresponding Levi factor. Let $\bfi_\bfP:\bfP\to\bfG$ and
$\bfa_\bfP:\bfP\to\bfM$ be the natural maps. Following Lusztig we define the restriction
functor $\Res^\bfG_\bfM:\calD(\bfG)\to\calD(\bfM)$ by setting
\eq{}
\Res^\bfG_\bfM(\calF)=(\bfa_\bfP)_!\bfi^*\calF
\end{equation}

%---------------------------------------------------------------------
\ssec{gsvolnoi}{The space $\tbG$}Let $\tbG$ denote the variety
of all pairs $(\bfB,g)$, where

$\bullet$\quad $\bfB$ is a Borel subgroup of $\bfG$,

$\bullet$\quad $g\in \bfB$.

One has natural maps $\alp: \tbG\to \bfT$ and $\pi:\ \tbG\to \bfG$
defined as follows. First of all, we set $\pi(\bfB,g)=g$. Now, in
order to define $\alp$, let us recall that for any Borel subgroup
$\bfB$ of $\bfG$ one has canonical identification $\mu_\bfB:\
\bfB/\bfU_{\bfB}\wt{\rightarrow} \bfT$, where $\bfU_{\bfB}$
denotes the unipotent radical of $\bfB$ (in fact, this is how the
abstract Cartan group $T$ is defined). Now we set
$\alp(\bfB,g)=\mu_\bfB(g)$.
%---------------------------------------------------------------------
\ssec{induction}{Induction}
Let now $\bfP=\bfB$. In this case $\bfM$ is the Cartan group $\bfT$.
We define the induction functor $\Ind_\bfT^\bfG:\calD(\bfT)\to\calD(\bfG)$
setting
\eq{}
\Ind_\bfT^\bfG(\calF)=\pi_!\alp^*[d](\frac{d}{2})
\end{equation}
where $\pi$ and $\alp$ are as above and $d=\dim\bfG-\dim\bfT$.

We set $\Spr=\Ind_\bfT^\bfG(\del_e)$ where $\del_e$ denotes the $\delta$-function
sheaf at the unit element of $\bfG$. It is known that $\Spr$ is
a perverse sheaf supported on the set $\calN$ of unipotent elements
of $\bfG$. Moreover, $\Spr$ is endowed with a natural $W$-action
(see for example \cite{BM}).

%---------------------------------------------------------------------------------------
\ssec{averaging}{A reformulation}For a subgroup $\bfH$ of $\bfG$ let $\calD^\bfH(\bfG)$ denote the derived
category of $\ell$-adic sheaves on $\bfG$ which are equivariant with respect
to the adjoint action.
Then the functor $\Ind_\bfT^\bfG$ can be rewritten as follows. Following
\cite{MV} let us define the averaging functor
$\avgb:\calD^\bfB(\bfG)\to\calD^\bfG(\bfG)$. Let $\del:\bfG\x \bfB\to\bfB$ and
$\eta:\bfG\x \bfB\to \bfG\underset \bfB\x \bfB$ be the natural maps
(in the definition of $\bfG\underset \bfB\x \bfB$ the group $\bfB$ acts by translations
on $\bfG$ and via the adjoint action on $\bfB$). Let also $\bfm_\bfB:
\bfG\underset \bfB\x\bfB\to\bfG$ be the map sneding a pair $(g,b)$ to $gbg^{-1}$ (note that
$\bfm_\bfB$ is proper). Let $\calF\in\calD^\bfB(\bfB)$. Then there exists canonical
$\calG\in\calD(\bfG\underset \bfB\x \bfG/\bfB)$ such that $\eta^*\calG=\del^*\calF$.
We define $\avgb(\calF)=(\bfm_\bfB)_!\calG$.

Abusing the notations we shall
denote the composition of $\avgb$  with the forgetful functor going from
$\calD^\bfG(\bfG)$ to $\calD(\bfG)$ by the same symbol.
Given any $\calF\in\calD(\bfT)$ its inverse image $\calF_\bfB$ to $\bfB$
with respect to the natural map $\bfB\to\bfT$ can be naturally regarded as
an object of $\calD^\bfB(\bfG)$. Then it is easy to see that
\eq{averaging}
\Ind^\bfG_\bfT(\calF)=\avgb(\calF_\bfB)[d](\frac{d}{2}).
\end{equation}
%---------------------------------------------------------------------------------------------------
We define now the functors $\tilInd_\bfT^\bfG:\calD(\bfT)\to\calD(\bfG)$ and
$\tilRes_\bfT^\bfG:\calD(\bfG)\to\calD(\bfT)$ by setting
\eq{}
\begin{align}
\tilInd_\bfT^\bfG(\calF)=\Ind_\bfT^\bfG\calF\ten H^*(\bfT,\qlb)\ten\qlbs^{\ten 2\dim \bfG}=\\
\Ind_\bfT^\bfG\calF\ten (H_c^*(\bfT,\qlb))^{\vee}[2d](d)
\end{align}
\end{equation}
where $(H_c^*(\bfT,\qlb))^{\vee}$ denotes the graded dual to $H_c^*(B)$
and
\eq{}
\tilRes_\bfT^\bfG\calG=\Res^\bfG_\bfT\calG\ten H^*_c(\bfG,\qlb)
\end{equation}

The following facts about the induction and restriction functors are essentially due to
Lusztig \cite{lu-char}, Theorem 4.4 and Ginzburg \cite{Gin}, Theorem 6.2. However, since these
results are stated in {\it loc. cit.} are stated only for character sheaves, we are going
to sketch the proofs.
%----------------------------------------------------------------------------------
\th{lus-gin}
\begin{enumerate}
\item
The functor $\Ind_\bfT^\bfG$ maps perverse sheaves to perverse sheaves.
\item
Let $\calG$ be a perverse sheaf on $\bfG$ which is equivariant with respect
to the adjoint action. Let also $\calF$ be any perverse sheaf on $\bfT$.
Then
\eq{adjointness}
\RHom(\tilRes^\bfG_\bfT\calG,\calF)=\RHom(\calG,\tilInd_\bfT^\bfG\calF)
\end{equation}
Moreover, for any $\calF,\calG$ as above the following diagram is commutative.
\begin{equation}
\begin{CD}
\RHom(\tilRes^\bfG_\bfT\calG,\calF)@>>>\RHom(\calG,\tilInd_\bfT^\bfG\calF)\\
@VVV                          @VVV\\
\RHom(\tilRes^\bfG_\bfT\fr^*\calG,\fr^*\calF)@>>>\RHom(\fr^*\calG,\tilInd_\bfT^\bfG\fr^*\calF)
\end{CD}
\end{equation}
(note that to write vertical arrows one needs to use the natural isomorphisms
$\tilInd_\bfT^\bfG\fr^*\calF\simeq \fr^*\tilInd_\bfT^\bfG\calF$ and
$\tilRes^\bfG_\bfT\fr^*\calG\simeq \fr^*\tilRes^\bfG_\bfT\calG$).
\item
Let $\calF$ be an irreducible perverse sheaf on $\bfT$. Assume that the
support of $\calF$ is a $W$-invariant subtorus in $\bfT$. Then for every $w\in W$ one has
a canonical isomorphism
\eq{3}
\Ind_\bfT^\bfG(\calF)\simeq \Ind_\bfT^\bfG(w^*\calF).
\end{equation}
\end{enumerate}
\eth

\noindent
{\it Remark.} One can show that \reft{lus-gin}(3) holds for any perverse sheaf
$\calF$ on $\bfT$. However, in this case the argument is a little more complicated
and we are not going to present it since we don't need it.
\prf
The first assertion of \reft{lus-gin} is standard (cf. \cite{lu-char}, Section 4.3).
Also the second assertion follows from standard adjointness properties of inverse and
direct images. Let us prove the third assertion.

Let $\bfT'\subset \bfT$ denote the support of $\calF$. Then we can find two reductive
groups ${\bfG_1}$, $\bfG_2$ and a surjective homomorphism $\kap:{\bfG_1}\x \bfG_2\to\bfG$
such that

1) the kernel of $\kap$ is a central subtorus in ${\bfG_1}\x\bfG_2$.

2) the preimage of $\bfT'$ under $\kap$ (with respect to some embedding of $\bfT$ into $\bfG$)
is equal to a maximal torus ${\bfT_1}$ in ${\bfG_1}$.

Let $\bfT_2$ be a (split) maximal torus in $\bfG_2$. We have the natural map
$\kap_\bfT:{\bfT_1}\x\bfT_2\to\bfT$ with connected kernel. It is easy to see that
it is enough to construct the isomorphism \refe{3} for $\kap_T^*\calF[\dim\ker \kap_\bfT]$.

Let $W_i$ denote the Weyl group of $\bfG_i$. Then $W=W_1\x W_2$.
On the other hand
$$
\Ind_{{\bfT_1}\x\bfT_2}^{{\bfG_1}\x\bfG_2}(\kap_T^*\calF[\dim\ker \kap_\bfT])=
\Ind_{\bfT_1}^{\bfG_1}(\kap_\bfT^*\calF[\dim\ker \kap_\bfT])
\boxtimes \Spr_{\bfG_2}.
$$
Since the second multiple is endowed with a natural action of $W_2$, it is enough to construct an
action of $W_1$ on the first multiple.
However, it follows from the fact that the map $\pi:\tbG\to\bfG$ is small that
$\Ind_{\bfT_1}^{\bfG_1}(\kap_T^*\calF[\dim\ker \kap_\bfT])$ is equal to the Goresky-MacPherson extension of
its restriction to the set of regular semi-simple elements in ${\bfG_1}$ where the construction
of the ${\bfG_1}$-equivariant structure is standard.
\epr

%-------------------------------------------------------------------------------------------
\ssec{}{Composition}Let $\bfP,\bfM$ be as above and let $W_\bfM\subset W$ be
the Weyl group of $\bfM$ (note that the embedding of $W_\bfM$ to $W$ depends
on $\bfP$).
\th{res-ind}
Let $\calF$ be an
irreducible perverse $W$-equivariant sheaf on $\bfT$ whose support
$\supp\, \calF$ is a $W$-stable subtorus in $\bfT$. Then
\eq{}
\Res^\bfG_\bfM\Ind_\bfT^\bfG(\calF)\simeq\Ind_{W_\bfM}^W(\Ind_\bfT^\bfM(\calF))
\end{equation}
and this isomorphism commutes with the natural actions of $W$ on both sides.
\eth

\noindent
{\it Remark.} For character sheaves this result is proved in \cite{lu-char}.

\prf
First of all consider $\Res^\bfG_\bfM\Ind_\bfT^\bfG(\calF)$ where $\calF$ is an arbitrary
object of $\calD(\bfT)$. We claim that it is glued from the complexes
$\Ind^\bfM_\bfT(w^*\calF)$ where $W$ runs over the representatives of the
cosets $W/W_\bfM$ of minimal length. Indeed, the sheaf $\Res^\bfG_\bfM\Ind_\bfT^\bfG(\calF)$
can be computed in the following way.

Consider the product $\bfZ=\bfP\underset{\bfG}\x\tbG$. We let $\del:\bfZ\to\bfT$ and
$\sig:\bfZ\to\bfM$ be the natural maps. Then $\Res^\bfG_\bfM\Ind_\bfT^\bfG(\calF)=
\sig_!\delta^*\calF[d](\frac{d}{2})$. On the other hand to every $w$ as
above there corresponds a locally closed stratum $\bfZ_w$ of $\bfZ$
(consisting of pairs $(\bfB,x\in\bfB\cap\bfP)$ where $\bfB$ and $\bfP$ are
in position $w\text{mod} W_\bfM$). We denote by $\del_w$ and $\sig_w$ the restrictions of
$\del$ and $\sig$ to $\bfZ_w$. Then looking at the Cousin complex associated with
the stratification $\bfZ_w$ we see that $\Res^\bfG_\bfM\Ind_\bfT^\bfG(\calF)$
is glued from the complexes $(\sig_w)_!\del_w^*\calF[d](\frac{d}{2})$.

Each $\bfZ_w$ has a natural map $\mu_w$ to ${\widetilde
\bfM}$. Namely, for every Borel subgroup $\bfB$ of $\bfG$ the image
of $\bfP\cap \bfB$ in $\bfM=\bfP/\bfU_\bfP$ is a Borel subgroup of $\bfM$
and we set $\mu_w(\bfB,x)=(\bfB\cap\bfP\mod\bfU_\bfP,x\mod\bfU_\bfP)$.
It is easy to see that $\mu_w$ is a locally trivial fibration with fiber
isomorphic to $\AA^{\dim\bfU_\bfP-l(w)}$. Also the composition of $\mu_w$
with the natural map $\alp_\bfM:{\widetilde \bfM}\to\bfT$ is equal to
$w\circ\delta$. This implies that
$(\sig_w)_!\del_w^*\calF[d](\frac{d}{2})\simeq \Ind_\bfT^\bfM(w^*\calF)$.
Hence $\Res^\bfG_\bfM\Ind_\bfT^\bfG(\calF)$  is glued from
the complexes $\Ind_\bfT^\bfM(w^*\calF)$.

Assume now that $\calF$ satisfies the conditions of the theorem.
Then $W$ acts naturally on $\Res^\bfG_\bfM\Ind_\bfT^\bfG(\calF)$ and it is easy to see that
it permutes the various subquotients $\Ind_\bfT^\bfM(w^*\calF)$ from which $\Res^\bfG_\bfM\Ind_\bfT^\bfG(\calF)$
is glued by the above argument. Hence $\Res^\bfG_\bfM\Ind_\bfT^\bfG(\calF)$ is isomorphic to
$\Ind_{W_\bfM}^W \Ind_\bfT$ .
\epr

Assume now that we are given two perverse sheaves $\calF_1$ and $\calF_2$ on $\bfT$
satisfying the conditions of \reft{res-ind}. Let also $\calG$ be a direct summand of
$\Ind_\bfT^\bfG(\calF_2)$. Then it follows from \reft{res-ind} that
$\Res^\bfG_\bfT(\calG)$ has a natural $W$-equivariant structure and
hence the same is true for $\tilRes^\bfG_\bfT(\calG)$. Therefore we have a natural
action of $W$ on $\RHom(\calF_1,\tilRes^\bfG_\bfT\calG)$. On the other hand, since
$W$ acts naturally on $\Ind_\bfT^\bfG\calF_1$ and on $H^*_c(\bfT,\qlb)$; hence it also acts on
$\tilInd_\bfT^\bfG\calF_1$. Therefore $W$ also acts on $\RHom(\tilInd_\bfT^\bfG\calF_1,\calG)$.
The proof of the following lemma is left to the reader.
%----------------------------------------------------------------------------------------------
\lem{commutes_with_W}
The isomorphism \refe{adjointness} commutes with the above $W$-actions.
\elem
%----------------------------------------------------------------------------------

%---------------------------------------------------------------------
Let $\bfU$ be the unipotent radical of a Borel subgroup $\bfB$ of $\bfG$.
Let $\bfr_\bfU:\bfG\to\bfG/\bfU$ denote the natural map.
\prop{ind-conv}
Let $\calG$ be a perverse sheaf on $\bfG$ which is  equivariant with respect to the adjoint action.

Assume that $(\bfr_\bfU)_!\calG$ is concentrated on $\bfT\subset \bfG/\bfU$.
Then for every $\calF\in\calD(\bfT)$ we have
\eq{conv-main}
\calG\star\Ind_\bfT^\bfG(\calF)\simeq \Ind_\bfT^\bfG(\Res^\bfG_\bfT\calG\star\calF)
\end{equation}
(in the above formula the restriction functor is taken with respect to the chosen
Borel subgroup $\bfB$).

The same result holds if we require that $(\bfr_\bfU)_*\calG$ vanishes outside
of $\bfT$ and replace $\star$-convolution by $*$-convolution.
\eprop
%--------------------------------------------------------------
\prf
We are going to prove the statement about $\star$-convolution. The proof
for $*$-convolution is analogous.

 It is easy to see that for any $\calH\in\calD^\bfB(\bfG)$
we have the natural isomorphism
\eq{conv-average}
\calG\star\avgb(\calH)[d](\frac{d}{2})\simeq \avgb(\calG\star\calH)[d](\frac{d}{2})
\end{equation}
Let us apply this to $\calH=\calF_\bfB$. Then the right hand side of
\refe{conv-average} is equal to the right
hand side of \refe{conv-main}.
On the other hand, the assumption that
$\bfr_\bfU\calG$ vanishes outside of $\bfT=\bfB/\bfU$ together with
$\bfU$-equivariance of $\calF_\bfB$ imply that
\eq{}
\calG\star\calF_\bfB=(\Res^\bfG_\bfT(\calG)\star\calF)_\bfB
\end{equation}
which finishes the proof.

\epr

%-----------------------------------------------------------------------------------------------
\sec{char}{Character sheaves and Deligne-Lusztig representations}
%-----------------------------------------------------------------------------------------------
\ssec{}{Maximal tori}Let us recall the classification of conjugacy classes
of $\fq$-rational maximal tori in $\bfG$. Recall that we assume that
$\bfG$ is split.

Let $\bfT$ denote the abstract Cartan group of $\bfG$ (with its canonical
$F$-rational structure).
Given $w\in W$ we can construct a new Frobenius morphism $\fr_w:\bfT\to
\bfT$ sending every $t\in\bfT$ to $w(\fr(t))$. In this way we get a new
$F$-rational structure on $\bfT$. We will denote the resulting torus by
$\bfT_w$. It is clear that if $w$ and $w'$ belong to the same
conjugacy class then $\bfT_w$ and $\bfT_{w'}$ are isomorphic.

The following result is proved in \cite{dl}:
\th{wind}
For every $w\in W$ there exists an embedding of $\bfT_w$ in $\bfG$ and
in this way we get a bijection between $G$-conjugacy classes of $F$-rational
maximal tori in $\bfG$ and conjugacy classes in $W$.
\eth
%----------------------------------------------------------
\ssec{}{Characters of tori}Let $\bfT$ be any algebraic torus over $F$ and
let $\theta$ be a character of $T$. One can associate to  $\theta$ an
$l$-adic local system $\calL_\theta$ on $\bfT$ in the following way.
Let $\alp:\bfT\to\bfT$ be the morphism given by
$\alp(t)=\fr(t)t^{-1}$. Then $\alp$ is a Galois \'etale
covering with Galois group equal to $T$. Hence
$T$ acts on the sheaf $\alp_*(\qlb)$. We set $\calL_\theta$ to be the
part of $\alp_*(\qlb)$ on which $T$ acts by means of $\theta$.

Let $w',w''\in W$ and let $\theta'$ (resp.
$\theta''$) be a character of $T_{w'}$ (resp. of $T_{w''}$). We say that
$\theta'$ and $\theta''$ are geometrically conjugate if there exists
$w\in W$ such that $w^*\calL_{\theta'}\simeq \calL_{\theta''}$
(note that as varieties over $\oF$ both $\bfT_{w'}$ and $\bfT_{w''}$
are identified with $\bfT$). This notion is introduced in \cite{dl}
in a slightly different language.
%---------------------------------------------------------------------------------------------
\ssec{}{Deligne-Lusztig representations}
Let $w\in W$ and let $\theta: T_w\to\qlb^\x$ be a character. In \cite{dl}
Deligne and Lusztig constructed a virtual representation $R_{\theta,w}$ of
$G$. We are going to need the following facts about $R_{\theta,w}$. Let
$\kg$ denote the Grothendieck group of representations of $G$. We have a
natural pairing $\<\, ,\>:\kg\ten\kg\to\ZZ$ such that if
$\pi_1,\pi_2\in\Irr(G)$ then
\eq{}
\< \pi_1,\pi_2\> =
\begin{cases}
1\qquad \text{if $\pi_1$ is isomorphic to $\pi_2$}\\
0\qquad \text{otherwise}
\end{cases}
\end{equation}

1) For every $\pi\in\Irr(G)$ there exists $w\in W$ and
$\theta:T_w\to\qlb^\x$ such that
\eq{}
\< \pi,R_{\theta,w}\>\neq 0
\end{equation}

2) One has
\eq{}
\< R_{\theta,w},R_{\theta',w'}\>\neq 0
\end{equation}
if and only if $\theta$ and $\theta'$ are geometrically conjugate.
%------------------------------------------------------------------------------------
\ssec{charsheaves}{Character sheaves}Let $\bfG$ be an arbitrary
reductive algebraic group over $F$.
Let us recall Lusztig's definition of
(some of) the character sheaves.

Let $\calL$ be a tame local system on $\bfT$. We define
$\chl=\Ind_\bfT^\bfG\calL$.
One knows (cf. \cite{lu-char}, \cite{laum})
that the sheaf $\chl$ is perverse.

\ssec{}{The Weil structure}
Assume now, that for some $w\in W$ there exists an
isomorphism $\calL\simeq \frw^*(\calL)$ and let us fix it.
It was observed by
G.~Lusztig in \cite{lu-char} that fixing such an isomorphism endows
$\chl$ canonically with a Weil structure (this follows immediately from
\reft{wind}).

Let now $\theta:T_w\to \qlb^{\x}$ be any character.  The following result is due to G.~Lusztig.
\th{lusztig}
\eq{}
\chi(\calK_{\calL_{\theta}})=q^{-\frac{d}{2}}\ch (R_{\theta,w})
\end{equation}
\eth
%----------------------------------------------------------------------------------------------------

%----------------------------------------------------------------------------------
\sec{torus}{$\gam$-sheaves on split tori}
%------------------------------------------------------------------------------
\ssec{}{}Let $\bfT$ be a split torus over $F$ and let
\eq{}
\rho=\lam_1\oplus...\oplus\lam_n
\end{equation}
be a good representation of $\td$. Recall that this means that there exists
a character $\sig:\bfT\to\GG_m$ such that
$\la\sig,\lam_i\ra>0$ for every $i$.

Each $\lam_i$ can be considered as a cocharacter of $\bfT$.
Let $\bfT_\rho=\GG_m^n$.
Define the map $\bfp_\rho:\bfT_\rho\to\bfT$ by setting
\eq{}
\bfp_\rho(t_1,...,t_n)=\lam_1(t_1)...\lam_n(t_n).
\end{equation}
We have
Let $\tr_\rho:\bfT_\rho\to \AA^1$ be given by
\eq{}
\tr_\rho(x_1,...,x_n)=x_1+...+x_n.
\end{equation}
Consider the complex $\phrot:=(\bfp_\rho)_!\tr_\rho^*\calL_\psi[n](\frac{n}{2})$ (on
$\bfT$).
%--------------------------------------------------------------------------
\th{torus-main}
Assume that $\lam_i$ is non-trivial for every $i=1,...,n$. Then
\begin{enumerate}
\item
The complex $\phrot$ is perverse.
\item
\eq{}
\text{supp}\,\phrot=\bfp_\rho(\bfT_\rho)
\end{equation}
\item
Assume in addition that $\rho$ is good.
Then the natural
map
\eq{shrik-star}
\phrot:=(\bfp_\rho)_!\tr_\rho^*\calL_\psi[n](\frac{n}{2})\to (\bfp_\rho)_*\tr_\rho^*\calL_\psi[n](\frac{n}{2})
\end{equation}
is an isomorphism and
$\phrot$ is an irreducible perverse sheaf on $\bfT$.
\end{enumerate}
\eth
\prf
Point (1) of \reft{torus-main} follows from the following
result.
%------------------------------------------------------------------------
\prop{conv-fourier}
Let $\lam:\GG_m\to\bfT$ be a non-trivial character and define
$\bPhi_\lam=\lam_*(\calL_\psi)[1](\frac{1}{2})$. Then the functors
$\bfC_{\lam,\star},\bfC_{\lam,*}: \calD(\bfT)\to\calD(\bfT)$ sending every $A\in\calD(\bfT)$ to
$\bPhi_\lam\star A$ and to $\bPhi_\lam * A$ respectively map
perverse complexes to perverse ones.
\eprop
%------------------------------------------------------------------------------
To see that \refp{conv-fourier} implies \reft{torus-main}(1) it is enough to
note that the complex $\bPhi_\lam$  is perverse for every $\lam$  and that
\eq{}
\phrot\simeq \bPhi_{\lam_1}\star ...\star\bPhi_{\lam_n}
\end{equation}

\noindent
{\it Proof of \refp{conv-fourier}.}
First of all there exists another torus $\bfT'_\lam$ together with an isogeny
$\bfq_\lam:\bfT'_\lam\to\bfT$ and an {\it injective} cocharacter $\lam':\GG_m\to \bfT$ such that
$\lam=\bfq_\lam\circ\lam'$. Let $\bfS_\lam=\bfT'_\lam/\GG_m$. Set also
$\obfT'_\lam=\bfT'_\lam\underset{\GG_m}\x \AA^1$.

We have the natural map $\bfs_\lam:\obfT'_\lam\to \bfS_\lam$ which endows
$\obfT'_\lam$ with a structure of a line bundle over $\bfS_\lam$.
Moreover, the dual vector bundle can be
naturally identified with $\bfs_{\lam^{-1}}:\obfT'_{\lam^{-1}}\to\bfS_{\lam^{-1}}=\bfS_\lam$. Thus one can consider the
Fourier-Deligne transform functor $\bfF_\lam:\calD(\obfT'_\lam)\to
\calD(\obfT'_{\lam^{-1}})$ (cf. e.g. \cite{KaLa}).

Let $j_\lam:\bfT'_\lam\to \obfT'_\lam$ denote the natural embedding.

The following lemma is straightforward from the definitions.
\lem{f-conv}
There exists a natural isomorphism of functors
\eq{shrik}
\bfj_{\lam}^*\circ\bfF_{\lam^{-1}}\circ (\bfj_{\lam^{-1}})_!\circ\circ \bfq_\lam^*
\simeq
\bfq_\lam^*\circ\bfC_{\lam,\star}
\end{equation}
and
\eq{star}
\bfj_{\lam}^*\circ\bfF_{\lam^{-1}}\circ (\bfj_{\lam^{-1}})_*\circ\circ \bfq_\lam^*
\simeq
\bfq_\lam^*\circ\bfC_{\lam,*}
\end{equation}
(note that we can take $\bfT'_\lam=\bfT'_{\lam^{-1}}$).
\elem
To see that \refl{f-conv} implies \refp{conv-fourier} it is enough to note
that a complex $A\in\calD(\bfT)$ is perverse if and only if $\bfq_\lam^*(A)$
is
perverse and that the functors $(\bfj_{\lam^{-1}})_!$ and $\bfF_{\lam^{-1}}$ map perverse
complexes to perverse ones (for the former this follows from the fact that
$\bfj_\lam$ is an affine open embedding and for the latter cf. \cite{KaLa}).
%----------------------------------------------------------------------------------
\ssec{}{Proof of \reft{torus-main}(2)}We are going to use induction on $n$. For $n=1$
the statement is obvious. Assume that we know the result for $n-1$.
Let $\rho'=\{ \lam_1,...,\lam_{n-1}\}$. Then we know that $\bPhi_{\bfT,\rho',\psi,!}$ is supported
on the image of $\bfp_\rho$. We have $\bPhi_{\bfT,\rho,\psi,!}=\bPhi_{\bfT,\rho',\psi,!}\star \bPhi_{\lam_n}$.
Thus our statement follows from \refl{f-conv} and from the following result.
%---------------------------------------------------------------------------------------------
\prop{support}
Let $\bfX$ be a scheme over $\oF$, $\pi:\bfL\to\bfX$ - a line bundle. Denote by $\bfL^{\vee}$
the dual line bundle and by $\bfj:\overset{\circ}\bfL\to\bfL$ the embedding of the complement to the
zero section (thus we have the natural isomorphism $\overset{\circ}\bfL\simeq \overset{\circ}\bfL^{\vee}$).
Let $\bfF_\psi:\calD(\bfL)\to\calD(\bfL^{\vee})$ denote the Fourier-Deligne transform corresponding
to the additive character $\psi$. Then for every $\calF\in\calD(\bfL)$ we have
\eq{}
\supp(\bfF_\psi(\bfj_!\calF))=\overline {\GG_m\cdot \supp(\calF)}
\end{equation}
(here bar denotes the Zariski closure).
\eprop
%---------------------------------------------------------------------------------------------
\prf The statement is immediately reduced to the case when $\bfX$ is a point.
Thus we have the embedding $\bfj:\GG_m\to \AA^1$ and a complex $\calF\in\calD(\GG_m)$.
We have to show that $\bfF_\psi(\bfj_!\calF)$ is non-zero at the generic point
of $\AA^1$ provided that $\calF\neq 0$.

Assume the contrary. Then the Grothendieck-Ogg-Schafarevich formula for the Euler-Poincar\'e characteristic
(cf. for example formula 2.3.1 in \cite{Katz}) implies that

1) $\calF$ is locally constant.

2) For any $t\in \oF^*$ the complex
$\calF\ten t^*\calL_\psi|_{\GG_m}$ is tame at infinity.

Clearly this is possible only if $\calF=0$.
\epr
%----------------------------------------------------------------------------------------

Let us now pass to the proof of \reft{torus-main}(3).
Let $\bfH$ denote the kernel of $\bfp_\rho$ and let $\del_\bfH$ be the constant sheaf
on $\bfH$ shifted by $\dim \bfH$ (regarded as a perverse sheaf on $\bfT_\rho$.

We must show that the natural map
$$
\phrot=(\bfp_\rho)_!\tr_\rho^*\calL_\psi[n](\frac{n}{2})\to (\bfp_\rho)_*\tr_\rho^*\calL_\psi[n](\frac{n}{2})
$$
is an isomorphism
and that $\phrot$ is an irreducible perverse sheaf. Taking inverse image to $\bfT_\rho$ we see that this
is equivalent to the following two statements:

1) The natural map
$\tr_\rho^*\calL_\psi[n](\frac{n}{2})\star \del_\bfH\to\tr_\rho^*\calL_\psi[n](\frac{n}{2})* \del_\bfH$
is an isomorphism.

2) The perverse sheaf $\tr_\rho^*\calL_\psi[n](\frac{n}{2})\star \del_\bfH[\dim \bfH]$ is irreducible
as an $\bfH$-equivariant
perverse sheaf (i.e. it has no $\bfH$-equivariant subsheaves).

Let $\bfj:\bfT_\rho=\GG_m^n\to\AA^n$ be the natural embedding. Let also $\bfF_{\AA^n}$ denote the Fourier-Deligne
transform on $\AA^n$. Then arguing as above it is easy to see that
\eq{}
\tr_\rho^*\calL_\psi[n](\frac{n}{2})\star \del_\bfH=\bfj^*\bfF_{\AA^n}(\bfj_!\del_\bfH)
\end{equation}
and
\eq{}
\tr_\rho^*\calL_\psi[n](\frac{n}{2}) *\del_\bfH=\bfj^*\bfF_{\AA^n}(\bfj_*\del_\bfH)
\end{equation}
Since $\bfF_{\AA^n}$ is an auto-equivalence which maps $\bfH$-equivariant perverse sheaves to
$\bfH$-equivariant perverse sheaves we see that in order to prove 1 and 2 above it is enough
to show that the map $\bfj_!\del_\bfH\to\bfj_*\del_\bfH$ is an isomorphism.

For this it is enough to show that $\bfH$ is closed in $\AA^n$. But $\bfH$ is a closed subset
of $\bfH_\sig=\ker \sig$. Hence it is enough to show that $\bfH_\sig$ is closed in $\AA^n$.
Let $a_i=\la \lam_i,\sig\ra,\ i=1,...,n$. Then
\eq{hsigma}
\bfH_\sig=\{ (t_1,...,t_n)\in \GG_m^n|\ t_1^{a_1}...t_n^{a_n}=1\}
\end{equation}
Since $\rho$ is good it follows that $a_i>0$ for every $i$. This together with \refe{hsigma}
clearly implies that $\bfH_\sig$ is closed in $\AA^n$.
\epr
%----------------------------------------------------------------------------------------------------
\ssec{}{Tame local systems}Recall that an $\ell$-adic local system $\calL$
on $\bfT$ is called {\it tame} if there exists a finite homomorphism
$\pi:\bfT'\to\bfT$ of some other algebraic torus $\bfT'$ to $\bfT$ such that
$\pi^*\calL$ is trivial.
%-------------------------------------------------------------------------------------
\th{torus-L}
Let $\rho=\oplus\lam_j$ be as above such that all $\lam_j$ are non-trivial.
Then for every tame local system $\calL$ on $\bfT$ one has
\eq{}
\dim H^i_c(\phrot\ten \calL)=
\begin{cases}
0\quad\text{if $i\neq 0$},\\
1\  \text{if $i=0$}.
\end{cases}
\end{equation}
The same is true for $H^i(\phrot\ten \calL)$.

We set $H_{\rho,\calL,\psi,!}:=H^0_c(\phrot\ten \calL^{-1})$ and
$H_{\rho,\calL,\psi,*}:=H^0(\phrot\ten \calL^{-1})$.
\eth
\prf
Let us prove the statement of \reft{torus-L} for
$H^i_c(\phrot\ten \calL)$. The proof for
$H^i(\phrot\ten \calL)$ is analogous.

First of all, we may assume without loss of generality that $\bfT=\bfT_\rho=\GG_m^n$,
and $\rho:\GG_m^n\to\gl(n)$ is the standard embedding. Indeed, the
definition of $\phrot$ and the projection formula imply that
the statements of \reft{torus-L} hold for $\calL$ if and only if they
hold for $\bfp_\rho^*\calL$.

On the other hand since $\calL$ is a one-dimensional tame local system on $\GG_m^n$ it follows
that there exist tame local systems $\calL_1, ...,\calL_n$ on $\GG_m$ such that
$\bfp_\rho^*\calL\simeq \calL_1\boxtimes ...\boxtimes \calL_n$.
Since in our case $\phrot=\calL_\psi\boxtimes ...\boxtimes \calL_\psi[n](\frac{n}{2})$,
\reft{torus-L}(1) follows from the following well-known lemma.
\lem{}
Let $\calL$ be a one-dimensional local system on $\GG_m$. Then
$H^i_c(\calL\ten\calL_\psi)=0$ if $i\neq 1$ and
$\dim H^1_c(\calL\ten\calL_\psi)=1$.
\elem

\epr
\cor{fignya}
Assume that $\rho$ is good. Then we have
\eq{}
\phrot\star\calL=H_{\rho,\calL,\psi,!}\ten \calL
\end{equation}
and
\eq{}
\phrot *\calL=H_{\rho,\calL,\psi,*}\ten \calL
\end{equation}
\ecor
%------------------------------------------------------------------------------------------
\sec{gln}{The basic example}
Now we assume again that $\bfG=\gln$.
Set $\phrog=\tr^*\calL_\psi[n^2](\frac{n^2}{2})$.
Let also $\phrot=\tr_\bfT^*\calL_\psi[n](\frac{n}{2})$ where
$\tr_\bfT:\bfT\to\AA^1$ is the restriction of the trace morphism
to the group $\bfT$ of diagonal matrices in $\gln$.
%-------------------------------------------------------------------
\th{convolution}
\begin{enumerate}
\item
One has an isomorphism of functors
\eq{four-com}
\phrog\star\Ind_\bfT^\bfG(\calF)\simeq \Ind_\bfT^\bfG(\phrot\star\calF)
\end{equation}
\item
Let $\calL$ be a tame local system on $\bfT$ endowed with an isomorphism
$\frw^*\calL\simeq\calL$ for some $w\in W$. In this case both sides of
\refe{four-com} are endowed with a natural Weil structure. Then the
isomorphism of \refe{four-com} is an isomorphism of Weil sheaves.
\end{enumerate}
\eth
%----------------------------------------------------------------------------------------
\prf
Let $\grg$ denote the Lie algebra of $\bfG$, i.e. the algebra of $n\x n$-matrices.
We have the natural embedding $\bfj_\bfG:\bfG\hookrightarrow \grg$.
Let us identify $\grg$ with its dual space by means of the form
$(x,y)\mapsto \tr(xy)$. Let $\bfF_\grg:\calD(\grg)\to\calD(\grg)$
denote the corresponding Fourier transform functor.
Then for every $\calG\in\calD(\bfG)$ we have
\eq{fourier_on_g}
\phrog\star\calG\simeq (\bfj_\bfG)^*\bfF_\grg((\bfj_\bfG)_!\calG^\iota)
\end{equation}
where  $\calG^\iota$ denotes the inverse image of $\calG$ with respect to
inversion map $g\mapsto g^{-1}$.
Similarly for any $\calF\in\calD(\bfT)$ we have
\eq{fourier_on_t}
\phrot\star\calF\simeq (\bfj_\bfT)^*\bfF_\grt((\bfj_\bfT)_!\calF^\iota)
\end{equation}

Let $\tilgrg$ be the space of all pairs $(\grb,x)$ where $\grb$ is a Borel subalagebra
in $\grg$ and $x\in\grb$. We have the natural open embedding $\bfj_\tbG:\tbG\hookrightarrow
\tilgrg$.

Let $\grt$ be the Cartan algebra of $\bfG$ (the Lie algebra of $\bfT$) and let
$\bfj_\bfT:\bfT\hookrightarrow\grt$ be the natural embedding.
Then  as in \refss{induction} we can define the induction
functor $\Ind_\grt^\grg:\calD(\grt)\to\calD(\grg)$. It follows immediately
from the definitions that for every $\calF\in\calD(\bfT)$ we have the natural isomorphisms
$$
(\bfj_\bfG)_!\Ind_\bfT^\bfG\calF\simeq \Ind_\grt^\grg((\bfj_\bfT)_!\calF)
$$
and
$$
\Ind_\bfT^\bfG \calF^\iota\simeq (\Ind_\bfT^\bfG\calF)^\iota.
$$
Hence the first statement of \reft{convolution} follows from the following
\lem{}
There is a natural isomorphism of functors
\eq{}
\bfF_\grg\circ\Ind_\grt^\grg\simeq \Ind_\grt^\grg\circ\bfF_\grt.
\end{equation}
\elem
\prf
Let $\pi:\tilgrg\to\grg$ and $\alp:\tilgrg\to\grt$ be the natural maps.
Then for every $\calF\in\calD(\grt)$ we have
$$
\Ind_\grt^\grg(\calF)=\pi_!\alp^*\calF[\dim\grg-\dim\grt](\frac{\dim\grg-\dim\grt}{2}).
$$
Let $\bfX$ denote the flag variety of $\bfG$. Then $\tilgrg$ can be regarded
naturally as a vector subbundle of the trivial vector bundle
$\bfX\x \grg$ over $\bfX$. Let $\eta:\bfX\x\grg\to\grg$ be the natural projection and
let also $\bfF_{\bfX\x\grg}$ denote the Fourier transform in the fiber of the vector
bundle $\bfX\x\grg$. It is known that in this situtation the functor $\eta_!$ commutes
with Fourier transform, i.e. there is a natural ismorphism of functors
$$
\eta_!\circ\bfF_{\bfX\x\grg}\simeq \bfF_\grg\circ\eta_!.
$$
Hence the lemma follows from the following observation:
for every $\calF\in\calD(\grt)$ we have
$\bfF_{\bfX\x\grg}(\alp^*\calF)\simeq \alp^*\bfF_\grt(\calF)$.
\epr
The second statement of \reft{convolution} is proved using similar considerations
and we leave it to the reader.
\epr
%--------------------------------------------------------------------------------------
\cor{}
Assume that an irreducible representation
$(\pi,V)$ appears in some Deligne-Lusztig representation $R_{\theta,w}$. Then
\eq{}
\gam_\psi(\pi)=\gam_{w,\psi}(\theta)
\end{equation}
where the notations are as in \refss{gln-int}.
In particular, $\gam(\pi)=\gam(\pi')$ if $\pi$ and $\pi'$
appear in the same virtual representation $R_{\theta,w}$.
\ecor
This follows immediately from \reft{convolution} and from \reft{lusztig}.
%-------------------------------------------------------------------------------------------
\sec{main}{$\gam$-sheaves: the main results}
%--------------------------------------------------------------------------------------

%---------------------------------------------------------------------------------------------------------------
\ssec{}{The perverse sheaf $\phrog$}
In what follows we assume that $\rho$ is good and $W$-equivariant.
We want to define a $W$-equivariant structure on $\phrot$.
Recall that we denote $\bfT_\rho=\GG_m^n$. Then the group
$W_\rho:=S_n$ acts naturally on $\bfT_\rho$.

Choose $w\in W$. We need to define an isomorphism
$\iota_w:w^*(\phrot){\wt\to}\phrot$. Let (as above) $w'$ be any
lift of $\rho(w)$ to $W_\rho$. Then one has
\eq{bred}
\bfp_\rho(w'(t))=w(\bfp_\rho(t)).
\end{equation}
The sheaf $\tr_\rho^*\calL_\psi$ is obviously
$W_\rho$-equivariant. This, together with \refe{bred}, gives rise
to an isomorphism $\iota_w':w^*(A_\rho){\wt\to}A_\rho$. We now
define
\eq{iota}
\iota_w:=(-1)^{l(w')-l(w)}\iota_w'.
\end{equation}
%----------------------------------------------------------------------------------------
\prop{equivariant}
\begin{enumerate}
\item
The isomorphism
$\iota_w$ does not depend on the choice of $w'$.
\item
The assignment
$w\mapsto\iota_w$ defines a $W$-equivariant structure on the sheaf
$\phrot$.
\end{enumerate}
\eprop
%-----------------------------------------------------------------------------------------
\prf
Clearly the second statement of \refp{equivariant} follows from the
first one. So, we just have to show that $\iota_w$ does not depend on the
choice of $w'$. For this it is enough to show the following.

$\bullet$ Let $s\in W_\rho$ be a simple reflection. Assume that
$\bfp_\rho\circ s=\bfp_\rho$ (i.e. $s$ is a lift of the unit element
$e\in W$ to $W_\rho$). Then $\iota_s$is equal to multiplication by
$-1$.

For this it is enough to prove the following lemma.
\lem{z2}
Consider the torus $\GG_m^2$ with coordinates $x$ and $y$.
Let $\bfm:\GG_m^2\to\GG_m$ be the multiplication map and let
$\bff:\GG_m^2\to\AA^1$ be given by $\bff(x,y)=x+y$. Then
the involution $s$ interchanging $x$ and $y$ acts on
$\bfm_!(\bff^*\calL_\psi)$ by means of multiplication by $-1$.
\elem
\prf
The quotient
of $\GG_m^2$ by the action of $\ZZ_2$ coming from $s$ is isomorphic to
$\AA^1\x \GG_m$ with coordinates $z=x+y$ and $t=xy$. Let $\bfq:\GG_m^2\to\GG_m^2/\ZZ_2$
be the natural map. Then it is enough to show that the direct image under
$\bfm$ of
$\bfq_!(\bff^*\calL_\psi)^{\ZZ_2}$ vanishes. But the latter sheaf is
isomorphic $\calL_\psi\boxtimes\qlb$ on $\GG_m^2/\ZZ_2\simeq \AA^1\x \GG_m$
and the required assertion follows from the fact that
$H^*_c(\AA^1,\calL_\psi)=0$.

\epr
\epr

Consider $\Ind_\bfT^\bfG(\phrot)$. It follows
from \reft{torus-main} that $\phrot$ satisfies the conditions
of \reft{wind}.
Define
\eq{}
\phrog=(\Ind_\bfT^\bfG(\phrot))^W
\end{equation}

\prop{}
\begin{enumerate}
\item
The sheaf $\phrog$ is non-zero and irreducible.
\item
Assume that $\im\bfp_\rho\cap\bfT_{rs}\neq \emptyset$ where $\bfT_{rs}$ denotes the set of
regular semi-simple elements in $\bfT$. Then $\phrog$ is equal to the
Goresky-MacPherson extension of its restriction to the set of regular
semisimple elements in $\bfG$. Moreover, the restriction of $\phrog$ to
the set $\bfG_r$ of regular elements in $\bfG$ is equal to
$\bfs^*(\bfq_!\phrot)^W[\dim\bfG-\dim\bfT](\frac{\dim\bfG-\dim\bfT}{2})$ where $\bfs:\bfG_r\to\bfT/W$ and
$\bfq:\bfT\to\bfT/W$ are the natural maps.
\end{enumerate}
\eprop
\prf
Let us first prove 2. The fact that $\phrog$ is equal to the
Goresky-MacPherson extension of its restriction to the set of regular
semisimple elements in $\bfG$ follows from smallness of the morphism
$\pi:\tbG\to\bfG$ (indeed, the fact that $\pi$ is small implies that
$\Ind_\bfT^\bfG\phrot$ is  equal to the
Goresky-MacPherson extension of its restriction to the set of regular
semisimple elements in $\bfG$ (since in this case
$\phrot$ is equal to the
Goresky-MacPherson extension of its restriction to the set of regular
semisimple elements in $\bfT$) and hence the same is true for any of its
direct summands).

Let us show that
$$
\phrog|_{\bfG_r}=\bfs^*(\bfq_!\phrot)^W[\dim\bfG-\dim\bfT](\frac{\dim\bfG-\dim\bfT}{2}).
$$
Both sides are equal to the Goresky-MacPherson extensions of their
restrictions to the set $\bfG_{rs}$ of regular semi-simple elements.
Hence it is enough to establish the above isomorphism on
$\bfG_{rs}$ where it is obvious.

Let us prove 1. Arguing as in the proof of \reft{lus-gin} we can assume that
one of the following holds:

(i) $\im \bfp_\rho=\bfT$.

(ii) $\im \bfp_\rho$ lies in $\bfZ(\bfG)$ (the center of $\bfG$).

Consider case (i). Then to show that $\phrog$ is irreducible it is enough
to show that $\phrog|_{\bfG_r}$ is irreducible (this follows from 2).
On the ohter hand since
$$
\phrog|_{\bfG_r}=\bfs^*(\bfq_!\phrot)^W[\dim\bfG-\dim\bfT](\frac{\dim\bfG-\dim\bfT}{2}).
$$
and since $\bfs$ has connected fibers it is enough to show that
$(\bfq_!\phrot)^W$ is irreducible which follows immediately from the
irreducibility of $\phrot$.

Let us now consider case (ii). In this case we claim that
$\phrog=\phrot$ (this means that both sheaves are supported on $\bfZ(\bfG)$ and are
equal there). Indeed, the sheaf $\Ind_\bfT^\bfG\phrot$ is supported on
$\bfZ(\bfG)\cdot\calN$ and  $\Ind_\bfT^\bfG\phrot=\phrot\boxtimes \Spr$.
Moreover, it follows from \refe{iota} that the action of $W$ on
$\Ind_\bfT^\bfG\phrot$ comes from the second multiple and there it is
equal to the standard action of $W$ on $\Spr$ twisted by the sign character.
It is well-known that the sheaf $\Hom_W(\text{sign},\Spr)$ is equal to
$\delta_e$ where $\delta_e$ denotes that $\delta$-function sheaf at the unit element
of $\bfG$. Hence we have
$$
\phrog=(\Ind_\bfT^\bfG\phrot)^W=\phrot\boxtimes (\Spr\ten \sign)^W=
\phrot\boxtimes\delta_e=\phrot.
$$
\epr

Let us now discuss the relation between the sheaf $\phrog$ and the function $\Phi_{G,\rho,\psi}$.
%------------------------------------------------------------------------------------------
\conj{main}
Let $\bfP\subset \bfG$ be a parabolic subgroup, $\bfU\subset \bfP$ -- its
unipotent radical, $\bfM=\bfP/\bfU$ -- the corresponding Levi group.
Let $\bfq_\bfP:\bfG\to\bfG/\bfU$ and $i_\bfP:\bfM\to \bfG/\bfU$ be the
natural morphisms. Assume that $\rho$ is good. Then
$(\bfq_\bfP)_!\phrog$ vanishes outside of $\bfM$.
\econj
%-----------------------------------------------------------------------------------------
\th{main}
 We have
$$
\Res^\bfG_\bfM \phrog=\bPhi_{\bfM,\rho,\psi}.
$$
Assume now that \refco{main}
holds. Then
\begin{enumerate}
\item
\eq{}
(\bfq_\bfP)_!\phrog\simeq (\bfi_\bfP)_! \bPhi_{\bfM,\rho,\psi}
\end{equation}
\item
Assume that $\rho_1$ and $\rho_2$ are good with respect to the same
character $\sig$ of $\bfG$. Then
\eq{}
\bPhi_{\bfG,\rho_1\oplus\rho_2,\psi}\simeq
\bPhi_{\bfG,\rho_1,\psi}\star\bPhi_{\bfG,\rho_2,\psi}\simeq
\bPhi_{\bfG,\rho_1,\psi}*\bPhi_{\bfG,\rho_2,\psi}
\end{equation}
\item
Let $\calL$ be a tame local system on $\bfT$. Then
\eq{}
\phrog\star \chl\simeq H_{\rho,\calL,\psi,!}\ten \chl
\end{equation}
and
\eq{}
\phrog * \chl\simeq H_{\rho,\calL,\psi,*}\ten \chl.
\end{equation}
If $\calL$ is endowed with an isomorphism $\frw^*\calL\simeq \calL$ then these
isomorphisms commute with the Weil structures on both sides
(note that due to $W$-equivariance of $\phrot$ every isomorphism $\frw^*\calL\simeq\calL$
endowes the spaces $H_{\rho,\calL,\psi,!}=H^0(\phrot\ten\calL^{-1})$  and
$H_{\rho,\calL,\psi,*}=H^0(\phrot\ten\calL^{-1})$ with a Frobenius action).
\end{enumerate}
\eth
\cor{sheaf-function}
Assume \refco{main} holds. Then
$$
\chi(\phrog)=\Phi_{G,\rho,\psi}.
$$
\ecor
\prf
Let $w\in W$, $\theta:T_w\to\qlb^\x$ and let $\calL=\calL_\theta$.
It follows from \reft{main}(3) and \reft{lusztig} that in order to prove \refc{sheaf-function}
it is enough to show that the scalar by which Frobenius acts on $H_{\rho,\calL,\psi,!}$ is equal
to $\gam_{w,\rho,\psi}(\theta^{-1})$. It follows from the definitions that it is enough
to do it in the case when $\bfG=\gl(n)$ and $\rho$ is the standard representation
where it is obvious.
\epr
%--------------------------------------------------------------------------------------------------------------
\conj{}
The functors $\calF\mapsto \calF\star \phrog$ and $\calF\mapsto
\calF*\phrog$ from $\calD(\bfG)$ to itself are exact with respect to the perverse
$t$-structure.
\econj
%---------------------------------------------------------------------------------------------
\th{main2}
\refco{main} holds for $\bfG$ of semi-simple rank $\leq 1$. In particular, \refc{sheaf-function}
holds for $\bfG$ of semi-simple rank $\leq 1$.
\eth
%------------------------------------------------------------------------------------------
\prf
Clearly we can assume that the semi-simple rank of $\bfG$ is equal to 1
(otherwise $\bfG$ is a torus and in this case there is nothing to prove).
Also, without loss of generality we may assume that  $\im \bfp_\rho\cap\bfT_{rs}\neq\emptyset$
(otherwise the support of $\phrog$ lies in the center of $\bfG$ and again
there is nothing to prove).

In this case the Weyl group $W$ is isomorphic to $\ZZ_2$. We denote by
$\sig$ the only non-trivial element in $W$. Let
also $\bfT'=\bfT/\bfT\cap [\bfG,\bfG]$ and let $\pi':\bfT\to\bfT'$
be the natural map. Thus $W\simeq \ZZ_2$ acts in the fibers of $\pi'$.
Hence we get a natural map $\pi:\bfT/W\to\bfT'$.

Let $\bfB$ be a Borel subgroup of $\bfG$ with unipotent radical $\bfU$.
Let $g\in\bfG$ such that $g\not\in\bfB$. Then $gu$ is a regular element of
$\bfG$ for every $u\in\bfU$. Moreover, the map $\bfs:\bfG_r\to\bfT/W$
identifies $g\bfU$ with one of the fibers of $\pi$.
Hence it is enough to show that $\pi_!\phrot=\pi_*\phrot=0$.
For this it enough to prve that $(\pi'_!\phrot)^W=(\pi'_*\phrot)^W=0$

First of all, it follows from \reft{torus-main} (applied to the torus
$\bfT'$) that $\pi_!\phrot=\pi_*\phrot$ and that the sheaf $\bPhi':=\pi'_!\phrot$
is irreducible. Hence $\sig\in W$ acts on $\bPhi'$ by means of multiplication
by a scalar. Since $\sig^2=1$ it follows that this scalar must be
$\pm 1$.

We claim that the above scalar is equal to $-1$. Since
$H^0(\bPhi')=H^0(\phrot)\neq 0$, it is enough to check that $\sig$
acts on $H^0(\phrot)$ by means of multiplication by $-1$.
Let $\sig'$ be a lift of $\sig$ to $S_n$ and let
$\iota'_\sig:\sig^*\phrot\simeq\phrot$ be the corresponding
isomorphism (we are using here the notations introduced before
\refp{equivariant}). Then it follows
from \refl{z2} that $\iota'_\sig$ induces multiplication by
$(-1)^{l(\sig')}$ on $H^0(\phrot)$. Hence
$\iota_\sig=(-1)^{l(\sig')-l(\sig)}\iota'_\sig=(-1)^{l(\sig')-1}\iota'_\sig$
acts on $H^0(\phrot)$ by means of multiplication by $-1$.
\epr
%---------------------------------------------------------------
\ssec{}{Proof of \reft{main}}
\ssec{proof-rest}{}Let us show that $\Res^\bfG_\bfM\phrog=\bPhi_{\bfM,\rho,\psi}$.
By \reft{res-ind} we have
$$
\Res^\bfG_\bfM\Ind_\bfT^\bfG(\phrot)=\Ind_{W_\bfM}^W\Ind_\bfT^\bfM(\phrot).
$$
By Frobenius reciprocity
$$
(\Ind_{W_\bfM}^W\Ind_\bfT^\bfM(\phrot))^W=(\Ind_\bfT^\bfM(\phrot))^{W_\bfM}=\bPhi_{\bfM,\rho,\psi}.
$$
Hence
$$
\Res^\bfG_\bfM\phrog=(\Res^\bfG_\bfM\Ind_\bfT^\bfG(\phrot))^W=\bPhi_{\bfM,\rho,\psi}.
$$
This clearly implies \reft{main}(1) if we assume that \refco{main} holds.

The fact that \refco{main} implies \reft{main}(3) is an immediate
consequence of \refp{ind-conv}.
Hence we just need to prove \reft{main}(2).
We will do that for $\star$-convolution. The proof for $*$-convolution is analogous.

Consider first $\bPhi_{\bfG,\rho_1,\psi}\star\Ind_\bfT^\bfG(\bPhi_{\bfT,\rho_2,\psi})$.
Then \refco{main} and \refp{ind-conv} imply that it is isomorphic
to
$$
\Ind_\bfT^\bfG(\Res^\bfG_\bfT\bPhi_{\bfG,\rho_1,\psi})\star\bPhi_{\bfT,\rho_2,\psi}\simeq
\Ind_\bfT^\bfG(\bPhi_{\bfT,\rho_1,\psi}\star\bPhi_{\bfT,\rho_2,\psi})\simeq
\Ind_\bfT^\bfG(\bPhi_{\bfT,\rho_1\oplus\rho_2,\psi}).
$$
Our assertion is obtained by taking $W$-invariants on both sides.
%-------------------------------------------------------------------------------
%----------------------------------------------------------------------------

%---------------------------------------------------------------------------------------

We conclude with the following
\th{gln-equality}
Let $\Phi_{G,\rho,\psi}'=\chi(\phrog)$. Then
\begin{enumerate}
\item
For every $w\in W$ and $\theta:T_w\to\qlb^\x$ the trace of $\Phi_{G,\rho,\psi}'$ in
$R_{\theta,w}$ is equal to that of $\Phi_{G,\rho,\psi}$.
\item
If $\bfG=\gl(n)$ then
\eq{gln-equality}
\Phi_{G,\rho,\psi}'=\Phi_{G,\rho,\psi}
\end{equation}
\end{enumerate}
\eth
\prf

Let us prove 1. First of all we know that
$\Tr(\Phi_{G,\rho,\psi},R_{\theta,w})=\gam_{w,\rho,\psi}\cdot\dim R_{\theta,w}$.
By \cite{dl}(Theorem 7.1) we have $\dim R_{\theta,w}=q^{-\frac{d}{2}}\frac{\# G}{\# T_w}$.
Hence we need to show that
\eq{idiotism}
\Tr(\Phi_{G,\rho,\psi}',R_{\theta,w})=\gam_{w,\rho,\psi}\cdot q^{-\frac{d}{2}}\frac{\# G}{\# T_w}.
\end{equation}

Let $\calL$ be the local system on $\bfT$ corresponding to
$\theta,w$. Thus $\calL$ is endowed with the natural isomorphism
$\fr_w^*\calL\simeq \calL$. We know that the character of $R_{\theta,w}$ is
equal to $\chi(\chl)q^{\frac{d}{2}}$. Also the inverse image of $\chl$ under the map
$g\mapsto g^{-1}$ is equal to $\calK_{\calL^{-1}}=\DD\calK_\calL$. Hence it follows
that $\Tr(\Phi_{G,\rho,\psi}',R_{\theta,w})$ is equal to
$$
q^{\frac{d}{2}}\sum (-1)^i \Tr(\fr,H^i_c(\phrog\ten\chl)).
$$
Recall that for any two compexes $A$ and $B$ on a scheme $\bfX$
we have
$$
(\Ext^*(A,B))^{\vee}=H^*_c(A\ten \DD B).
$$
Hence it follows from
\refe{adjointness} that
we have a canonical isomorphism
$$
H^*_c(\phrog\ten\DD \tilInd_\bfT^\bfG\calL)\simeq H^*_c(\tilRes\phrog\ten\calL^{-1}).
$$
Changing $\calL$ to $\calL^{-1}$ and taking into account the canonical isomorphism
$\DD\Ind_\bfT^\bfG\calL^{-1}\simeq \Ind_\bfT^\bfG\calL$ we see that we have a canonical
isomorphism
$$
H^*_c(\phrog\ten\Ind_\bfT^\bfG\calL)[-2d](-d)\ten H^*_c(\bfT,\qlb)\simeq
H^*_c(\Res\phrog\ten\calL)\ten H^*_c(\bfG,\qlb)
$$
In other words, since we know that $\Res^\bfG_\bfT\phrog=\phrot$,
we get the isomorphism
\eq{idiot}
H^*_c(\phrog\ten\Ind_\bfT^\bfG\calL)[-2d](-d)\ten H^*_c(\bfT,\qlb)\simeq
H^*_c(\phrot\ten\calL)\ten H^*_c(\bfG,\qlb).
\end{equation}
It is easy to see that this isomorphism commutes with the action of Frobenius
on both sides where the actions of Frobenius on $H^*_c(\phrog\ten\Ind_\bfT^\bfG\calL)$
and $H^*_c(\bfG,\qlb)$ are standard and the actions on $H^*_c(\bfT,\qlb)$ and
$H^*_c(\phrot\ten\calL)$ are via $\fr_w$.
Thus \refe{idiotism} follows from \refe{idiot} by taking traces of Frobenius.

Let us prove the second assertion.
It is known that if $\bf=\gl(n)$ then the virtual Deligne-Lusztig representations $R_{\theta,w}$
generate over $\QQ$ the Grothendieck group of the category of finite-dimensional
representations of $G$  (cf. for example the introduction to
\cite{lu-book}). Hence \refe{gln-equality} follows from
\reft{gln-equality}(1).
\epr

\end{document}